\documentclass[12pt]{article}
\usepackage{amsmath,amsthm}
\usepackage{amssymb,latexsym}
\usepackage{enumerate}
\usepackage[english]{babel}
\usepackage{graphicx}

\def\Q{{\mathbb Q}}
\def\Z{{\mathbb Z}}

\newtheorem{lemma}{Lemma}
\newtheorem{theorem}[lemma]{Theorem}

\title{
Binomial Thue equations and \\ power integral bases in pure quartic fields
}
\author{
Istv\'{a}n Ga\'{a}l\thanks{
        This work was partially supported by the European Union 
        and the European Social Fund through project Supercomputer, 
        the national virtual lab (grant no.: TAMOP-4.2.2.C-11/1/KONV-2012-0010)
        and also supported in part by K67580 and K75566 from the
        Hungarian National Foundation for Scientific Research,
         }\;
and L\'aszl\'o Remete
\\
University of Debrecen, Mathematical Institute \\
            H--4010 Debrecen Pf.12., Hungary \\
            e--mail: igaal@science.unideb.hu, hagyma42@gmail.com \\ \\
}

\begin{document}

\maketitle
\thispagestyle{empty}

\renewcommand{\thefootnote}{}

\footnote{2010 \emph{Mathematics Subject Classification}: Primar 11R16; Secondary 11D59, 11Y50}

\footnote{\emph{Key words and phrases}: pure quartic fields, power integral bases, supercomputers}

\renewcommand{\thefootnote}{\arabic{footnote}}
\setcounter{footnote}{0}

\begin{abstract}
It is a classical problem in algebraic number theory to decide if a number
field admits power integral bases and further to calculate all generators of 
power integral bases. This problem is especially delicate to consider in
an infinite parametric family of number fields. In the present paper we
investigate power integral bases in the infinite parametric family of
pure quartic fields $K=\Q({}^4\sqrt{m})$. 

We often pointed out close connection of various types of Thue equations
with calculating power integral bases (cf. \cite{thueapple}, \cite{book}). 
In this paper we determine power integral bases in
pure quartic fields by reducing the corresponding index form equation to
quartic binomial Thue equations. We apply recent results on the solutions of
binomial Thue equations of type $a^4-mb^4=\pm 1\; ({\rm in}\; a,b\in\Z)$
as well as we perform an extensive calculation by a high performance computer
to determine "small" solutions (with $\max(|a|,|b|)<10^{500}$) 
of binomial Thue equations for $m<10^7$. 

Using these results on binomial Thue equations we characterize power integral bases
in infinite subfamilies of pure quartic fields. 
For $1< m< 10^7$ we give all generators of power integral bases with coefficients
$<10^{1000}$ in the integral basis.
\end{abstract}

\section{Introduction}

Let $K$ be an algebraic number field of degree $n$ 
with ring of integers $\Z_K$. 
There is an extensive literature (cf. \cite{book}) of {\em power integral bases}
in $K$, that is integral bases of the form $\{ 1,\gamma,\gamma^2,\ldots,\gamma^{n-1}\}$.
In general, if $\{1,\omega_2,\ldots,\omega_n\}$ is an integral basis of $K$, then
the discriminant of the linear form $\omega_2 X_2+\ldots+\omega_n X_n$
can be written as
\[
D_{K/Q}(\omega_2 X_2+\ldots+\omega_n X_n)=(I(X_2,\ldots,X_n))^2 D_{K}
\]
where $D_K$ denotes the discriminant of the field $K$,
and $I(X_2,\ldots,X_n)$ is the {\em index form} 
corresponding to the above integral basis.
As is known, for any primitive element $\gamma=x_1+\omega_2 x_2 +\ldots +\omega_n x_n\in\Z_K$
(that is $K=\Q(\gamma)$) we have
\[
I(\gamma)=(\Z_K^+:\Z[\gamma]^+)=|I(x_2,\ldots,x_n)|,
\]
where the index of the additive groups of the corresponding rings are taken.
$I(\gamma)$ is the {\em index of} $\gamma$, which does not depend on $x_1$.
Therefore the element $\gamma=x_1+\omega_2 x_2 +\ldots +\omega_n x_n$
generates a power integral basis of $K$ if and only if
$x_1\in\Z$ and $(x_2,\ldots,x_n)\in \Z^{n-1}$ is a solution of the
{\em index form equation} 
\begin{equation}
I(x_2,\ldots,x_n)=\pm 1 \;\;\; {\rm with}\;\;\; x_2,\ldots,x_n\in\Z .
\label{ifeq}
\end{equation}
This is the way the problem of power integral bases reduces to the
resolution of diophantine equations. The elements $\beta,\gamma\in\Z_K$
are called {\it equivalent}, if $\beta\pm \gamma\in\Z$. In this case
$I(\beta)=I(\gamma)$. We only consider generators of power integral bases
up to equivalence.

There is an extensive literature of index form equations and power integral bases.
Until today there are efficient algorithms for calculating generators
of power integral bases in number fields of degree up to 6, as well as
for some higher degree fields of special structure. (The interested reader
should consult I.Ga\'al \cite{book} for a general overview, for fields up to degree 5 
and for special higher degree fields. For sextic fields see 
Y.Bilu, I.Ga\'al and K.Gy\H ory \cite{degreesix}.)
Considering this problem not only in a single number field but in an infinite family
of fields is an especially challenging task. However we have succeeded
in several cases, see e.g. I.Ga\'al \cite{complex4}, I.Ga\'al and M.Pohst \cite{c5}, 
I.Ga\'al and T.Szab\'o \cite{gsz1}, \cite{gsz2}, I.Ga\'al \cite{book}.

\section{Pure quartic fields}

The purpose of the present paper is to consider power integral bases in pure quartic fields. 
Recently pure cubic fields were studied by I.Ga\'al and T. Szab\'o \cite{purecubic}.
Therefore it was straithforward to closely consider the problem of
power integral bases in pure quartic fields that has not been investigated before.
The main goal of this paper is to show that the problem of power integral bases 
in pure quartic fields can be reduced to quartic binomial Thue equations.

Let $m$ be a positive integer and consider the pure quartic field $K=\Q({}^4\sqrt{m})$.
Set $\alpha={}^4\sqrt{m}$. Assume that $m$ is square free, $m\neq\pm 1$ and
$m\equiv 2,3 \;(\bmod 4)$. According to 
A.Hameed, T.Nakahara, S.M.Husnine and S.Ahmad \cite{hnha}
$\{1,\alpha, \alpha^2, \alpha^3\}$ is an integer basis in $K$.
We shall use this integer basis to construct the index form equation
therefore in the following we continue to assume these conditins on $m$.

Using results on binomial Thue equations and making direct calculations
we obtained the following statements.

\vspace{0.5cm}

As we shall see it is easy to deal with negative values of $m$.

\begin{theorem}
Let $m<-1$ be a square free integer with $m\equiv 2,3 \; (\bmod 4)$.
Then up to equivalence $\alpha={}^4\sqrt{m}$ is the only generator of 
power integral bases in $K=\Q({}^4\sqrt{m})$.
\label{th1}
\end{theorem}

Next we consider values of $m$ composed of at most two primes less than
or equal to 29.

\begin{theorem}
Let $m>1$ be a square free integer with $m\equiv 2,3 \; (\bmod 4)$
composed of at most two primes $p,q$ with $2\leq p,q\leq 29$.
Set $\alpha={}^4\sqrt{m}$ and let $K=\Q({}^4\sqrt{m})$.
Then up to equivalence the only generators of 
power integral bases in $K=\Q({}^4\sqrt{m})$ are 
$\vartheta=x\alpha+y\alpha^2+z\alpha^3$ with\\
$(x,y,z)=(1,0,0)$ for arbitrary $m$,\\
$(x,y,z)=(4,\pm 2,1)$ for $m=15$,\\
$(x,y,z)=(25,\pm 10,4)$ for $m=39$.\\
\label{theorem13}
\end{theorem}

We construct an infinite series of values of $m$ for which
the corresponding binomial Thue equation is easy to solve therefore
we can characterize the power integral bases in $K=\Q({}^4\sqrt{m})$.

\begin{theorem}
Let $s,t$ be positive integer parameters and let
\[
m=\frac{(s^4t\pm 1)^4-1}{s^4}
\]
be a square free integer with $m>1$ and $m\equiv 2,3 \; (\bmod 4)$.
Set $\alpha={}^4\sqrt{m}$.
Then up to equivalence the only generators of 
power integral bases in $K=\Q({}^4\sqrt{m})$ are 
$\vartheta=x\alpha+y\alpha^2+z\alpha^3$ with\\
$(x,y,z)=(1,0,0), ((s^4t\pm 1)^2, \pm(s^4t\pm 1) s, s^2)$.
\label{abab}
\end{theorem}

The following theorem uses results on binomial Thue
equations for small values of $m$.

\begin{theorem}
Let $1<m< 2000$ be a square free integer, with $m\equiv 2,3 \; (\bmod 4)$.
For $400<m< 2000$ we also assume that $m$ is odd.
Set $\alpha={}^4\sqrt{m}$ and let $K=\Q({}^4\sqrt{m})$.
Then up to equivalence the only generators of 
power integral bases in $K=\Q({}^4\sqrt{m})$ are 
$\vartheta=x\alpha+y\alpha^2+z\alpha^3$ with \\
$(x,y,z)=(1,0,0)$ for arbitrary $m$,\\
$(x,y,z)=(1,\pm 1,1)$ for $m=2$,\\
$(x,y,z)=(4,\pm 2,1)$ for $m=15$,\\
$(x,y,z)=(25,\pm 10,4)$ for $m=39$,\\
$(x,y,z)=(9,\pm 3,1)$ for $m=82$,\\
$(x,y,z)=(16,\pm 4,1)$ for $m=255$,\\
$(x,y,z)=(121,\pm 22, 4)$ for $m=915$,\\
$(x,y,z)=(36,\pm 6,1)$, for $m=1295$.\\
\label{b1theorem}
\end{theorem}

Our last theorem was obtained by an extensive calculation using a supercomputer
of "small" solutions of binomial Thue equations.
The range $m<10^7$ and the bound $10^{1000}$ for the solutions 
must be sufficient for practical applications.

\begin{theorem}
Let $1<m< 10^7$ be a square free integer with $m\equiv 2,3 \; (\bmod 4)$.
Set $\alpha={}^4\sqrt{m}$ and let $K=\Q({}^4\sqrt{m})$.
Then up to equivalence the only generators $\vartheta=x\alpha+y\alpha^2+z\alpha^3$
of power integral bases in $K=\Q({}^4\sqrt{m})$ with
\[
\max(|x|,|y|,|z|)<10^{1000}
\]
are $(x,y,z)=(1,0,0)$ for arbitrary $m$, moreover
the following $(x,y,z)$ for the given values of $m$:
\[
\begin{array}{|c|c|c|c|}
\hline
 m&x&\pm y&z\\ \hline
 2& 1& 1 & 1\\ \hline
 15& 4& 2  &  1\\ \hline
 39& 25 & 10 & 4 \\ \hline
 82& 9 & 3 & 1 \\ \hline
 255& 16& 4 & 1\\ \hline
 410& 81& 18& 4  \\ \hline
 626& 25& 5 & 1\\ \hline
 915& 121&22&4\\ \hline
 1295& 36& 6& 1\\ \hline
 2402& 49&7& 1 \\ \hline
 6562& 81& 9& 1 \\ \hline
 12155& 441& 42& 4\\ \hline
  14642& 121& 11& 1\\ \hline
 17490& 529&46&4\\ \hline
 20735& 144&12& 1\\ \hline
 24414& 625&50&4\\ \hline
 28562& 169&13& 1\\ \hline
 33215& 729& 54&4 \\ \hline
38415& 196& 14& 1\\ \hline
 50626& 225&15& 1\\ \hline
 61535& 3969& 252&16\\ \hline
 65535& 256& 16& 1\\ \hline
 83522& 289&17& 1\\ \hline
 130322& 361&19& 1\\ \hline
 144590& 1521&  78&4\\ \hline
 159999& 400& 20& 1\\ \hline
 194482& 441&21& 1\\ \hline
 234255& 484&22& 1\\ \hline
 279842& 529&23& 1\\ \hline
 390626& 625&25& 1\\ \hline
 505679& 6400& 240& 9\\ \hline
 531442& 729&27& 1\\ \hline
    \end{array}
\hspace{1cm}
\begin{array}{|c|c|c|c|}
\hline
 m&x&\pm y&z\\ \hline
  707282& 841&29& 1\\ \hline
 757335& 3481& 118& 4\\ \hline
  809999& 900&30& 1\\ \hline
 923522& 961&31& 1\\ \hline
   1081730 & 16641& 516&16\\ \hline
1185922 & 1089 & 33 & 1\\ \hline
1336335 & 1156 & 34 & 1\\ \hline
1416695 & 4761 & 138 & 4\\ \hline
1500626 & 1225 & 35 & 1\\ \hline
1679615 & 1296 & 36 & 1\\ \hline
1755519 & 33124 & 910 & 25\\ \hline
1874162 & 1369 & 37 & 1\\ \hline
1977539 & 5625 & 150 & 4\\ \hline
2313442 & 1521 & 39 & 1\\ \hline
2559999 & 1600 & 40 & 1\\ \hline
2825762 & 1681 & 41 & 1\\ \hline
3111695 & 1764 & 42 & 1\\ \hline
3262539 & 7225 & 170 & 4\\ \hline
3418802 & 1849 & 43 & 1\\ \hline
3580610 & 7569 & 174 & 4\\ \hline
3851367 & 196249 & 4430 & 100\\ \hline
4100626 & 2025 & 45 & 1\\ \hline
4879682 & 2209 & 47 & 1\\ \hline
5764802 & 2401 & 49 & 1\\ \hline
6765202 & 2601 & 51 & 1\\ \hline
7034430 & 10609 & 206 & 4\\ \hline
7311615 & 2704 & 52 & 1\\ \hline
7596914 & 11025 & 210 & 4\\ \hline
7890482 & 2809 & 53 & 1\\ \hline
8503055 & 2916 & 54 & 1\\ \hline
8608519 &105625 & 1950 & 36\\ \hline
9150626 & 3025 & 55 & 1\\ \hline
9834495 & 3136 & 56 & 1\\ \hline
  \end{array}
\]
\label{hpc}
\end{theorem}

\section{Auxiliary tools I: Index form equations in quartic fields}

Calculating generators of power integral bases is equivalent
to solving the index form equation (\ref{ifeq}). For quartic number fields
it is an equation of degree 6 in 3 variables. However by using the method of
I.Ga\'{a}l, A.Peth\H{o} and M.Pohst \cite{gppsim}
(see also I.Ga\'{a}l \cite{book}) this can be reduced to cubic and
quartic Thue equations that can be solved by using standard technics (cf. \cite{book}).
We cite here this result in detail since this will be a basic tool in our proofs.

Let $K$ be a quartic number field generated by a root $\xi$ with minimal polynomial 
$f(x)=x^4+a_1x^3+a_2x^2+a_3x+a_4\in\Z[x]$.
We represent any $\vartheta\in\Z_K$ in the form
\begin{equation}
\vartheta=\frac{1}{d}\left( a+x\xi+y\xi^2+z\xi^3  \right)
\label{alpha8}
\end{equation}
with coefficients $a,x,y,z\in\Z$ and with a common denominator $d\in\Z$.
Let $n=I(\xi)$,
\[
F(u,v)=u^3-a_2u^2v+(a_1a_3-4a_4)uv^2+(4a_2a_4-a_3^2-a_1^2a_4)v^3
\]
a binary cubic form over $\Z$ and
\begin{eqnarray*}
Q_1(x,y,z)&=&x^2-xya_1+y^2a_2+xz(a_1^2-2a_2)+yz(a_3-a_1a_2)+z^2(-a_1a_3+a_2^2+a_4)\\
Q_2(x,y,z)&=&y^2-xz-a_1yz+z^2a_2
\end{eqnarray*}
ternary quadratic forms over $\Z$.

\begin{lemma} (I.Ga\'{a}l, A.Peth\H{o} and M.Pohst \cite{gppsim})\\
If $\vartheta$ of (\ref{alpha8}) generates a power integral basis in $K$
then there is a solution $(u,v)\in \Z^2$ of
\begin{equation}
F(u,v)=\pm \frac{d^{6}m}{n}
\label{res}
\end{equation}
such that
\begin{eqnarray}
u&=&Q_1(x,y,z), \nonumber \\
v&=&Q_2(x,y,z).   \label{Q12}
\end{eqnarray}
\label{lemma4}
\end{lemma}

\section{Auxiliary tools II: Binomial Thue equations}

Let $g>0$ be a given integer. 
The index form equation in pure quartic fields
will be reduced to a {\em binomial Thue
equation} of type
\begin{equation}
a^4-gb^4=\pm 1 \;\; ({\rm in}\; a,b\in\Z).
\label{bthue}
\end{equation}
In order to solve the index form equation we shall use some corresponding 
results on binomial Thue equations.

If we know a positive solution of the binomial Thue equation, then
the following result of M.A.Bennett is useful to exclude any other positive solutions.

\begin{lemma}(M.A.Bennett \cite{b2001})\\
If $f,g$ are integers with $fg\neq 0$ and $n\geq 3$ then the equation
\[
fa^n-gb^n=\pm 1 \;\;\; ({\rm in}\; a,b\in\Z)
\]
has at most one solution in positive integers $(a,b)$.
\label{bennettlemma}
\end{lemma}

The following theorem describes all solutions of binomial Thue equation if its coefficients
are composed of at most two primes less than or equal to 29.

\begin{lemma} (K.Gy\H ory and \'A.Pint\'er \cite{b2012}) \\
If $f<g$ are positive integers composed of at most two primes $p,q$ 
with $2\leq p,q\leq 29$ and $n\geq 3$ then the only solutions of the equation
\[
fa^n-gb^n=\pm 1  \;\; ({\rm in}\;\; a,b\in\Z)
\]
are those with 
\[
f\in\{1,2,3,4,7,8,16\}, a=b=1,
\]
and 
\[
n=3, (f,a)=(1,2),(1,3),(1,4),(1,8),(1,9),(1,18),(1,19),(1,23),
\]
\[
(2,2),(3,2),(5,3),(5,11),(11,6),
\]
\[
n=4, (f,a)=(1,2),(1,3),(1,5),(3,2),
\]
\[
n=5, (f,a)=(1,2),(1,3),(9,2),
\]
\[
n=6, (f,a)=1,2).
\]
\label{lemma13}
\end{lemma}

The following results of A.Bazs\'o, A.B\'erczes, K.Gy\H ory and \'A.Pint\'er \cite{b2010}
describe the solutions of binomial Thue equations with $f=1$.
For simplicity we quote the case $n=4$ only, that is we consider equation
(\ref{bthue}).

\begin{lemma}(A.Bazs\'o, A.B\'erczes, K.Gy\H ory and \'A.Pint\'er \cite{b2010})\\
If $1<g\leq 400$ then all positive integer solutions $(a,b)$, with $|ab|>1$ of
equation (\ref{bthue}) are
\[
(g,a,b)=(5,3,2),(15,2,1),(17,2,1),(39,5,2),
\]
\[
(80,3,1),(82,3,1),(150,7,2),(255,4,1),(257,4,1). 
\]
\label{b1lemma}
\end{lemma}

\begin{lemma}(A.Bazs\'o, A.B\'erczes, K.Gy\H ory and \'A.Pint\'er \cite{b2010})\\
If $400<g\leq 2000$ is odd, then all positive integer solutions $(a,b)$, with $|ab|>1$ of
equation (\ref{bthue}) are
\[
(g,a,b)=(915,11,2), (1295,6,1), (1297,6,1), (1785,13,2).
\]
\label{b2lemma}
\end{lemma}

Remark. The solutions with $m=82,915,1295$ of Lemmas \ref{b1lemma}, \ref{b2lemma}
are indicated in the Corregindum \cite{b2010corr} to \cite{b2010}.

Our last theorem was obtained by an extensive computation
using a high performance computer. We calculate the "small" 
solutions up to $10^{500}$ of binomial Thue equations (\ref{bthue}) for $g$ up to $10^7$. 
Note that in several cases the calculation of fundamental units of the field 
$K=\Q({}^4\sqrt{g})$ is extremely time consuming. 
For example computing the fundamental unit in $K=\Q({}^4\sqrt{18071})$ 
takes more than 10 minutes using Kash \cite{kant} and the coefficients of the 
fundamental unit in the integral basis have ca. 1000 digits.
Therefore using the present
machinery it does not seem feasible to completely solve all these Thue equations.
However, according to our experience, these solutions yield all solutions with very high probability.

\begin{theorem}
All solutions with 
\[
\max(|a|,|b|)<10^{500}
\]
of the binomial Thue equation (\ref{bthue}) for $1<g<10^7$ 
(with irreducible left hand side) are listed in the following table.
\label{hpcbthue}
\end{theorem}

{\scriptsize

\[
\begin{array}{|c|c|c|}
\hline
 g&a& b \\ \hline 
        2  &   1  &   1  \\  \hline    
        5  &   3  &   2  \\  \hline        
        15  &   2  &   1  \\  \hline        
        17  &   2  &   1  \\  \hline        
       39  &   5  &   2  \\  \hline          
        80  &   3  &   1  \\  \hline          
        82  &   3  &   1  \\  \hline           
      150  &   7  &   2  \\  \hline            
       255  &   4  &   1  \\  \hline           
       257  &   4  &   1  \\  \hline             
      410  &   9  &   2  \\  \hline               
       624  &   5  &   1  \\  \hline              
       626  &   5  &   1  \\  \hline              
     915  &   11  &   2  \\  \hline              
      1295  &   6  &   1  \\  \hline              
      1297  &   6  &   1  \\  \hline               
     1785  &   13  &   2  \\  \hline              
      2400  &   7  &   1  \\  \hline               
      2402  &   7  &   1  \\  \hline                
     3164  &   15  &   2  \\  \hline                 
      4095  &   8  &   1  \\  \hline                  
      4097  &   8  &   1  \\  \hline                   
     5220  &   17  &   2  \\  \hline                   
      6560  &   9  &   1  \\  \hline                     
      6562  &   9  &   1  \\  \hline                     
     7140  &   239  &   26  \\  \hline                 
     8145  &   19  &   2  \\  \hline                       
     9999  &   10  &   1  \\  \hline                      
    10001  &   10  &   1  \\  \hline                   
    12155  &   21  &   2  \\  \hline                  
    14640  &   11  &   1  \\  \hline                
    14642  &   11  &   1  \\  \hline                 
    17490  &   23  &   2  \\  \hline                   
    20735  &   12  &   1  \\  \hline                 
    20737  &   12  &   1  \\  \hline                 
    24414  &   25  &   2  \\  \hline                    
    28560  &   13  &   1  \\  \hline                  
    28562  &   13  &   1  \\  \hline                  
    33215  &   27  &   2  \\  \hline                 
    38415  &   14  &   1  \\  \hline                   
    38417  &   14  &   1  \\  \hline                    
    44205  &   29  &   2  \\  \hline                   
    50624  &   15  &   1  \\  \hline                   
    50626  &   15  &   1  \\  \hline                    
    57720  &   31  &   2  \\  \hline                     
   61535  &   63  &   4  \\  \hline                           
\end{array}  
\hspace{0.5cm}
\begin{array}{|c|c|c|}
\hline
 g&a& b \\ \hline 
    65535  &   16  &   1  \\  \hline                     
    65537  &   16  &   1  \\  \hline                     
   69729  &   65  &   4  \\  \hline                    
    74120  &   33  &   2  \\  \hline     
    83520  &   17  &   1  \\  \hline                     
    83522  &   17  &   1  \\  \hline                    
    93789  &   35  &   2  \\  \hline                   
   104975  &   18  &   1  \\  \hline                   
    104977  &   18  &   1  \\  \hline                   
114240  &   239  &   13  \\  \hline               
   117135  &   37  &   2  \\  \hline                   
    130320  &   19  &   1  \\  \hline                 
    130322  &   19  &   1  \\  \hline                
   144590  &   39  &   2  \\  \hline                  
    159999  &   20  &   1  \\  \hline              
    160001  &   20  &   1  \\  \hline                    
   176610  &   41  &   2  \\  \hline                     
    194480  &   21  &   1  \\  \hline                  
    194482  &   21  &   1  \\  \hline                  
   213675  &   43  &   2  \\  \hline                   
    234255  &   22  &   1  \\  \hline                   
    234257  &   22  &   1  \\  \hline                    
   256289  &   45  &   2  \\  \hline                  
    279840  &   23  &   1  \\  \hline                     
    279842  &   23  &   1  \\  \hline                    
   304980  &   47  &   2  \\  \hline                    
    331775  &   24  &   1  \\  \hline                      
    331777  &   24  &   1  \\  \hline                 
   360300  &   49  &   2  \\  \hline                    
    390624  &   25  &   1  \\  \hline                  
    390626  &   25  &   1  \\  \hline                    
  422825  &   51  &   2  \\  \hline                      
    456975  &   26  &   1  \\  \hline                    
    456977  &   26  &   1  \\  \hline                    
   493155  &   53  &   2  \\  \hline                  
   505679  &   80  &   3  \\  \hline                    
518440  &   161  &   6  \\  \hline                   
    531440  &   27  &   1  \\  \hline                   
    531442  &   27  &   1  \\  \hline                     
544685  &   163  &   6  \\  \hline                       
   558175  &   82  &   3  \\  \hline                     
   571914  &   55  &   2  \\  \hline                     
    614655  &   28  &   1  \\  \hline                  
    614657  &   28  &   1  \\  \hline                   
   659750  &   57  &   2  \\  \hline                   
   707280  &   29  &   1  \\  \hline                      
\end{array}  
\hspace{0.5cm}
\begin{array}{|c|c|c|}
\hline
 g&a& b \\ \hline 
     707282  &   29  &   1  \\  \hline                     
   757335  &   59  &   2  \\  \hline                       
    809999  &   30  &   1  \\  \hline                     
    810001  &   30  &   1  \\  \hline                      
   865365  &   61  &   2  \\  \hline                     
    923520  &   31  &   1  \\  \hline                     
    923522  &   31  &   1  \\  \hline    
     984560  &   63  &   2  \\  \hline                       
   1016190  &   127  &   4  \\  \hline                    
   1048575  &   32  &   1  \\  \hline                    
   1048577  &   32  &   1  \\  \hline                   
1081730  &   129  &   4  \\  \hline                    
1115664  &   65  &   2  \\  \hline                       
   1185920  &   33  &   1  \\  \hline                       
   1185922  &   33  &   1  \\  \hline                   
  1259445  &   67  &   2  \\  \hline                       
   1336335  &   34  &   1  \\  \hline                        
   1336337  &   34  &   1  \\  \hline                       
  1416695  &   69  &   2  \\  \hline                   
   1500624  &   35  &   1  \\  \hline                      
   1500626  &   35  &   1  \\  \hline                        
  1588230  &   71  &   2  \\  \hline                         
  1679615  &   36  &   1  \\  \hline                        
   1679617  &   36  &   1  \\  \hline                  
1755519  &   182  &   5  \\  \hline                       
  1774890  &   73  &   2  \\  \hline                      
   1874160  &   37  &   1  \\  \hline                       
   1874162  &   37  &   1  \\  \hline                        
  1977539  &   75  &   2  \\  \hline                         
  2085135  &   38  &   1  \\  \hline                     
   2085137  &   38  &   1  \\  \hline                       
2197065  &   77  &   2  \\  \hline                           
  2313440  &   39  &   1  \\  \hline                       
   2313442  &   39  &   1  \\  \hline                      
  2434380  &   79  &   2  \\  \hline                     
   2559999  &   40  &   1  \\  \hline                       
   2560001  &   40  &   1  \\  \hline                     
2690420  &   81  &   2  \\  \hline                      
  2825760  &   41  &   1  \\  \hline                    
   2825762  &   41  &   1  \\  \hline                 
2966145  &   83  &   2  \\  \hline                     
  3111695  &   42  &   1  \\  \hline                          
   3111697  &   42  &   1  \\  \hline                          
  3262539  &   85  &   2  \\  \hline                        
   3418800  &   43  &   1  \\  \hline                        
   3418802  &   43  &   1  \\  \hline                                    
\end{array}  
\hspace{0.5cm}
\begin{array}{|c|c|c|}
\hline
 g&a& b \\ \hline 
  3580610  &   87  &   2  \\  \hline      
  3748095  &   44  &   1  \\  \hline                     
   3748097  &   44  &   1  \\  \hline                    
3851367  &   443  &   10  \\  \hline                    
3921390  &   89  &   2  \\  \hline                         
   4100624  &   45  &   1  \\  \hline                       
   4100626  &   45  &   1  \\  \hline                   
4285935  &   91  &   2  \\  \hline                        
  4477455  &   46  &   1  \\  \hline                      
   4477457  &   46  &   1  \\  \hline                    
  4675325  &   93  &   2  \\  \hline   
     4879680  &   47  &   1  \\  \hline                     
   4879682  &   47  &   1  \\  \hline                         
  5090664  &   95  &   2  \\  \hline                      
5198685  &   191  &   4  \\  \hline                   
   5308415  &   48  &   1  \\  \hline                   
   5308417  &   48  &   1  \\  \hline                     
5419875  &   193  &   4  \\  \hline                     
  5533080  &   97  &   2  \\  \hline                       
   5764800  &   49  &   1  \\  \hline                     
   5764802  &   49  &   1  \\  \hline                    
6003725  &   99  &   2  \\  \hline                       
   6249999  &   50  &   1  \\  \hline                  
   6250001  &   50  &   1  \\  \hline                   
6503775  &   101  &   2  \\  \hline                    
  6765200  &   51  &   1  \\  \hline                     
   6765202  &   51  &   1  \\  \hline                    
7034430  &   103  &   2  \\  \hline                      
    7311615  &   52  &   1  \\  \hline                   
  7311617  &   52  &   1  \\  \hline                      
7596914  &   105  &   2  \\  \hline                   
   7890480  &   53  &   1  \\  \hline                     
   7890482  &   53  &   1  \\  \hline                   
8192475  &   107  &   2  \\  \hline                   
8295040  &   161  &   3  \\  \hline                    
8398565  &   323  &   6  \\  \hline                 
   8503055  &   54  &   1  \\  \hline                    
   8503057  &   54  &   1  \\  \hline                    
8608519  &   325  &   6  \\  \hline                     
8714960  &   163  &   3  \\  \hline                      
8822385  &   109  &   2  \\  \hline                       
   9150624  &   55  &   1  \\  \hline                     
   9150626  &   55  &   1  \\  \hline                     
9487940  &   111  &   2  \\  \hline                    
  9834495  &   56  &   1  \\  \hline                       
   9834497  &   56  &   1  \\  \hline                      
\end{array}
\]

}

\section{Proofs}

Let $\alpha={}^4\sqrt{m}$ ($m$ is square free, $m\neq\pm 1$, with $m\equiv 2,3 \; (\bmod 4)$).
Then by the result of A.Hameed, T.Nakahara, S.M.Husnine and S.Ahmad \cite{hnha} 
$\{1,\alpha, \alpha^2, \alpha^3\}$ is an integer basis 
in the pure quartic field $K=\Q(\alpha)$. This means $n=I(\alpha)=1$ and any $\vartheta\in\Z_K$
can be represented in the form (\ref{alpha8}) with $d=1$. Therefore the equations
(\ref{res}) and (\ref{Q12}) can be written in the following form:
\begin{eqnarray}
u(u^2+4mv^2)&=&\pm 1,   \label{uv}\\
x^2-mz^2&=&u,           \label{q1}\\
y^2-xz&=&v.             \label{q2}
\end{eqnarray}
Equation (\ref{uv}) implies $u=\pm 1$ and $u^2+4mv^2=\pm 1$, whence $v=0$.

{\bf Case I.} Let $m<0$.
Equation (\ref{q1}) implies $x=1$ and $z=0$. By equation (\ref{q2}) we get $y=0$.
This yields that up to equivalence $\vartheta=\alpha$ is the only generator of power 
integral bases. This proves Theorem \ref{th1}.

{\bf Case II.} Let $m>1$. 
Equation (\ref{q1}), that is $x^2-mz^2=\pm 1$ implies that $\gcd(x,z)=1$.
Equation (\ref{q2}), that is $y^2=xz$ implies that $x$ and $z$ are both squares.
(Note that if $(x,y,z)$ is a solution then so also is $(-x,-y,-z)$, therefore
we may assume that $x,z$ are both nonnegative.)

Set $x=a^2$ and $z=b^2$ with $a,b\in\Z$. Then equation (\ref{q1}) yields the 
binomial Thue equation 
\begin{equation}
a^4-mb^4=\pm 1 \;\; ({\rm in}\; a,b\in\Z).
\label{binom}
\end{equation}
All solutions (up to sign) of (\ref{q1}) and (\ref{q2}) are obtained by
\begin{equation}
x=a^2,\;\; y=\pm ab\;\; z=b^2
\label{xyzab}
\end{equation}
where $a,b$ are solutions of (\ref{binom}).

\vspace{0.5cm}

Theorem \ref{theorem13} is a consequence of Lemma \ref{lemma13}.

\vspace{0.5cm}

To prove Theorem \ref{abab} observe that 
$(a,b)=(s^4t\pm 1,s)$ is a positive solution of equation (\ref{binom}).
By Lemma \ref{bennettlemma} this and $(a,b)=(1,0)$ 
are all solutions of this equation.
$(x,y,z)$ are calculated from $(a,b)$ by (\ref{xyzab}).
\vspace{0.5cm}

Theorem \ref{b1theorem} is a consequence of Lemma \ref{b1lemma} and of Lemma \ref{b2lemma}.

\vspace{0.5cm}

Theorem \ref{hpcbthue} was proved by calculating "small" solutions
of binomial Thue equations (\ref{binom}) with
$\max(|a|,|b|)<10^{500}$ by using the method of Peth\H o \cite{pet}. 
The procedure was implemented in Maple \cite{maple}, we used 1200 digits accuracy. 
This extensive calculation involving almost $10^7$ binomial
Thue equations was performed on the supercomputer (high performance computer) 
network situated in Debrecen-Budapest-P\'ecs-Szeged in Hungary. 
The total running time was about 120 hours calculated for a single node which yields 
a few hours using parallel computing with a couple of nodes. 

Theorem \ref{hpc} is a consequence of Theorem \ref{hpcbthue} using (\ref{xyzab}).

\end{document}